\documentstyle[12pt]{article}
\newcommand{\C}{1\!\!\!C}

\newcommand{\R}{I\!\!R}
\newcommand{\Sc}{{\cal S}}
\newcommand{\be}{\begin{eqnarray*}}
\newcommand{\ee}{\end{eqnarray*}}

      \date{}
\begin{document}

\setcounter{page}{1}


\title
{
The  $:\!\phi^4_4\!:$ quantum field theory, I.
\\ Wave operator, holomorphity and Wick kernel.
\footnote{This work
 is supported partly by RBRF 96-01649.
It is the second paper of the
project $\phi^4_4\cap M.$
}
}
\author{Edward  P. Osipov
\\  Department of Theoretical Physics
\\ Sobolev Institute of Mathematics
\\  630090 Novosibirsk, RUSSIA
\thanks{E-mail address: osipov@math.nsk.su}
}
      \date{}
\maketitle
\begin{flushright}

funct-an/9605002
\end{flushright}

\medskip

\medskip
\begin{abstract}
With the help of the complex  structure and the wave
operator of the nonlinear classical Klein-\-Gordon
equation with the interaction $u^4_4$ we define the
Wick kernel of the interacting quantum field in
four-\-dimensional space-\-time and consider its
properties. In particular, the diagonal of
this Wick kernel is (real) solutions of the
classical nonlinear
 Klein-\-Gordon
equation with the interaction $u^4_4.$
\end{abstract}

\section*{1. Introduction.}

In the paper
\cite{Osi94a},
see also
\cite{Hei74,Rac75},
we have  constructed
 the solution\footnote{
The constructed  solution is nontrivial
\cite[Theorem 4.3]{Osi94a},
\cite{Osi96c},
but here
we do not discuss the question of the nontriviality of the
$:\!\phi^4_4\!:$ quantum theory
\cite{Fro82,Aiz82,AizG83,Cal88,Hei74,PedSZ92,Osi94a}.
}
of the quantum wave equation
 $$
\Box \phi(t,x) + m^2\phi(t,x) +
\lambda:\!\phi^3(t,x)\!:=0\eqno(1.1)
$$
in four--dimensional space--time.
In the integral form this equation is the
following
 $$
\phi(t,x)=\phi_{in}(t,x)-\lambda\int^t_{-\infty}
\int\;R(t-\tau,x-y):\!\phi^3(\tau,y)\!:d\tau d^3 y,
\eqno(1.2)
$$
$$
\phi_{out}(t,x)
= \phi_{in}(t,x) - \lambda\int^{+\infty}_{-\infty}
\int\;R(t-\tau,x-y):\!\phi^3(\tau,y)\!:d\tau d^3 y,
\eqno(1.3)
$$
where the $in$-field $\phi_{in}$ is the $in$-coming
free quantum field and $:\,\,\,:$ is the Wick
normal ordering with respect to the free
field $\phi_{in}.$

In the papers
\cite{Hei74,Osi94a}
 this solution has been correctly defined as
a bilinear form on
$D_{coh}(\vartheta)\times D_{coh}(\vartheta),$ where
$D_{coh}(\vartheta)$ is a subspace in the Fock space of
the $in$--field $\phi_{in}$.
The subspace $D_{coh}(\vartheta)$ is generated by linear
combinations of coherent vectors which are near to the vacuum.

In the present paper we consider
the solution of the quantum nonlinear wave equation
(1.1),
which was constructed in
\cite{Hei74,Osi94a},
 from the other point of view.

We introduce explicitly the
Wick kernel for the
 interacting field
$$
\phi(e_{z_1},e_{z_2})
={1\over 2}\exp(\langle z_1,
z_2\rangle_{H^{1/2}(\R^3,\C)})
 (RW_{in}R^{-1}(\bar z_1+z_2)
+ \overline{RW_{in}R^{-1}
(\overline{\bar z_1+z_2})} )
\eqno(1.4)
$$
and for the $out$-going field
$$
\phi_{out}(e_{z_1},e_{z_2})
={1\over 2}\exp(\langle z_1,
z_2\rangle_{H^{1/2}(\R^3,\C)})
 (RSR^{-1}(\bar z_1+z_2)
+ \overline{RSR^{-1}
(\overline{\bar z_1+z_2})} ),
\eqno(1.5)
$$
where $R$ is the isomorphism,
$$
R:\quad H^1(\R^3)\oplus L_2(\R^3) \to
 H^{1}(\R^3,\C),
$$
and $W_{in}, W_{out}$ are the $in$- and $out$-wave
operator, respectively, of the classical
nonlinear equation with the interaction of the
4th degree. This Wick symbol is defined for all
 coherent vectors with
finite energy.
Due to holomorphity
the  expression
(1.4)
is holomorphic
 in $\bar z_1, z_2$
for $z_1,z_2\in H^1(\R^3,\C),$
 see
\cite{Osi95a}.

The expression (1.4)
 gives,  in fact,  the mathematical
realization of a physical idea
to write an interacting field
with the help of a wave operator.
 Due to holomorphity     (1.4) has the
  quite correct sense written as
$$
\phi={1\over 2}\left(:RW_{in}R^{-1}(\phi_{in}):
+ :RW_{out}R^{-1}(\phi_{in}):\right),
$$
i.e.
$$
\phi = \sum^\infty_{n=1} \phi_n(:\phi_{in}...\phi_{in}:)
 = \sum^\infty_{n=1} {1\over 2}
(R_n + \overline{R}_n)(:\phi_{in}...\phi_{in}:),
\eqno(1.6)
$$ where $R_n$ is the unique
generalized function from $\Sc(\R^{3n},\C),$
which corresponds to the $n$-th Frech{\'e}t
derivative of the function
$RW_{in}R^{-1}$ at zero in the
complex Hilbert space  $H^1(\R^3,\C)$
and
 $\overline{R}_n$ is its complex conjugation,
i.e. (we denote later $W_{in}$ as
 $W$ without the subscript $``in"$)
$$
R_n(f_1\otimes...\otimes f_n) =
{1\over n!}d^n
RWR^{-1}(0)(f_1,...,f_n), $$
$f_j\in\Sc(\R^{3n},\C).$  In other words,
$R_n$ is equal to the
$n$-th
derivative   of the function
 $RWR^{-1}$ at zero
and the derivative is taken over
the positive frequency part
of $in$-solution.
The positive frequency part
\footnote{
We use the usual ``physical" terminology for this notion.
It is defined uniquely for the massive linear
Klein--Gordon equation, i.e. for free solutions.
}
of $in$-solution may be defined with the help of the
isomorphism $R$ which gives the one-to-one
correspondence with the initial $in$-data.

The expansion (1.6) is defined by
  the Taylor expansion of the
(classical) wave operator at zero and it is
convergent
on coherent vectors near to the vacuum.
The expansion (1.6) defines the
bilinear form and
the holomorphity allows to prove that the Wick
kernel of (1.6) can be extended uniquely
 on all coherent vectors with
finite energy (and is equal to (1.4)\,).

Thus, we construct the
$:\!\phi^4_4\!:$ quantum field theory
along the following pathway.
To construct the
 $:\!\phi^4_4\!:$ quantum field theory
 we define the Wick kernel for the interacting
 quantum field as (1.4).
As a consequence we obtain that the interacting
 quantum field satisfies the expansion (1.5)
in the sense of bilinear forms
on an appropriate subspace of the Fock space.
The complex structure and estimates allow us to
prove that (1.4) is defined correctly as a Wick
kernel in the sense of
Paneitz, Pedersen, Segal, and Zhou
\cite{PanPSZ91}.

Separately, we prove that the
solution of Eq. (1.3), defined as a
 bilinear form, has the Wick kernel
equal to (1.4) on
coherent vectors near to the vacuum, i.e. in
a zero neighbourhood,
\cite{Osi95b},
see also
\cite{Osi94a,Hei74,Rac75}.
 In addition, this Wick kernel (1.4)
defines the unique
 operator--valued generalized function from
 $\Sc^1 (\R^4) \times \Sc^1 (\R^4)$ into
$L(\Phi o\kappa).$
Here
 $\Sc^\alpha (\R^4)$ are the Gelfand
spaces of test functions
\cite[v.~2]{GelS58},
\cite[v.~4]{GelV61},
 $L(\Phi o\kappa)$
is the space of linear bounded operators
 in the space
   $\Phi o\kappa$
  and $\Phi o\kappa$
denotes the Fock space
of the free quantum $in$-field.
Moreover, this operator-valued generalized function
can be extended on a Gelfand space
 $\Sc^{\alpha},$ $\alpha < 6/5,$
and these Gelfand spaces with $\alpha >1$ contain
dense subspaces of functions with compact supports in
coordinate space
\cite{Osi96b}.
 This result allows us to construct Wightman functions
(as generalized function of Jaffe type)
and matrix elements of the quantum scattering
operator.
It is possible to consider  Wightman axioms:
 the conditions of positivity,
 spectrality, Poin\-car{\'e} invariance, locality,
 cluster property, and
 asymptotic
completeness
\cite{Osi96c}.

In this paper we introduce the Wick kernel
(1.4) for the interacting quantum field of
the $:\!\phi^4_4\!:$ quantum field theory and
consider its properties.



\section*{2. Complex structure.}

In this section we formulate and prove the assertions
connected with the complex structure of solutions
of the classical (real) nonlinear wave
equation with the cubic
nonlinearity in four--dimensional Minkowski space--time,
$$
 \Box u+m^2 u+\lambda u^3=0,\;\; m>0,\;\;
\lambda>0. \eqno(2.1) $$
We  use these statements for the
construction of the $:\phi^4_4:$ quantum field theory.
 To introduce the Wick kernel of the interacting
quantum field of the
 $:\!\phi^4_4\!:$ theory
 we need the assertion about holomorphity
(Theorem 2.1)
and about
 $T$ (and/or
$CPT$) symmetry (Theorem 2.2).

The solution of Eq. (2.1) is uniquely determined by its
Cauchy data
$$ u |_{t=0} = \varphi, \quad\partial_t u|_{t=0} = \pi.
$$
Let $W_{in}$ $(= W)$ be the incoming
wave operator (the $in$--wave
operator, or the backward wave operator) of Eq. (2.1).
The operator $W_{in}$  maps $in$-data
at time zero into initial data at time zero. Let
$W_{out}$
 be the outgoing wave operator
(the $out$-wave operator, or the forward wave operator).
The operator $W_{out}$ maps  $out$--data into
 initial data.
Let
 $S$  be  the scattering operator of the classical nonlinear
wave equation (2.1).

\medskip
{\bf Theorem 2.1}
\cite[Theorem 1.1]{Osi95a}.

{\it Let $R(\varphi,\pi)
= \varphi + i \mu^{-1} \pi,$ $\mu=(-\Delta + m^2)^{1/2},$
be the invertible map of
$\Sc_{Re}(\R^3) \oplus \Sc_{Re}(\R^3)$
onto $\Sc(\R^3,\C).$ This map $R$ defines also
isomorphisms $(H^{1/2}(\R^3) \oplus H^{-1/2}(\R^3), J)$)
onto
$H^{1/2}(\R^3,\C)$ and of $(H^1(\R^3) \oplus L_2(\R^3), J)$
onto
$H^1(\R^3,\C).$
Here
$$ J=R^{-1}iR =
\pmatrix{0&-\mu^{-1}\cr\mu &0\cr}.
$$

The map $RWR^{-1}$
is defined correctly as the mapping  from
 $H^1(\R^3,\C)$ onto
 $H^1(\R^3,\C)$ and is
the complex analytic  mapping of the complex
Hilbert space
 $H^1(\R^3,\C)$ onto itself.
In particular, for
$z_{in}(\alpha)
= \sum^N_{j=1}\alpha_j z_{in,j},$
$\alpha_j \in \C,$
$z_{in,j} \in H^1(\R^3,\C)$,
$h \in H^{1/2}(\R^3,\C),$ the functions
$\langle h,RWR^{-1}
(z_{in}(\alpha))\rangle_{H^{1/2}(\R^3,\C)}$
are entire holomorphic function on
$(\alpha_1,...\alpha_N) \in \C^N.$

The same assertions are valid for the transformation
$RW_{out}R^{-1}$ and $RSR^{-1}.$
}

\smallskip
\medskip
Theorem 2.1 is proved in the paper
\cite[Theorem 1.1]{Osi95a}.

\medskip
\medskip
Let  $\cal P$ be the Poincar\'e group.
Let  ${\cal P}_0$ be the connected component of the
Poincar\'e group.
Given $(a,\Lambda)\in {\cal P}$ and a finite--energy
solution $u$ of (2.1), define
$(a,\Lambda)u$ by
$$
((a,\Lambda)u)(t,x) = u((a,\Lambda)^{-1}(t,x)).
$$
Then
$(a,\Lambda)u$ is also a finite--energy solution
of (2.1), and the map
$((a,\Lambda), u)\to (a,\Lambda)u$ defines an
action of the group ${\cal P}$ on the Hilbert
space
$H^1\oplus L_2$ of initial data.
We denote this map by $U(a,\Lambda).$ When
the coupling constant in (2.1) $\lambda=0,$ i.e.
for the case of the linear Klein--Gordon equation, the maps
 $U_0(a,\Lambda)$ are linear.

The wave operators have the following basic properties.

\medskip
{\bf Theorem 2.2.}

{\it
The wave operator $W_{in}$ and $W_{out}$ are
correctly defined on the Hilbert space
$H^1\oplus L_2.$
Moreover, they intertwine the free and interacting action
of the
Poincar\'e group, i.e.
$$
U(a,\Lambda)=W_{in}U_0(a,\Lambda)W_{in}^{-1},
$$
 $$
U(a,\Lambda)=W_{out}U_0(a,\Lambda)W_{out}^{-1}
$$
for all $(a,\Lambda) \in {\cal P}_0.$
}

\medskip
{\it Proof of Theorem  2.2.}
See, for instance,
\cite{Bae92}.
Theorem  2.2 is proved.

\medskip
\medskip
Let $\Theta^T$ be the time reflection operator on
the space of initial data,
$\Theta^T(\varphi,\pi) = (\varphi,-\pi).$
Let $\Theta^P$ be the space reflection operator on
the space of initial data,
 $\Theta^P(\varphi,\pi) = (\varphi^P,\pi^P),$
where $\varphi^P(x) = \varphi(-x),$
 $\pi^P(x)=\pi^P(-x).$

\medskip
{\bf Theorem 2.3 ($T$-symmetry).}

{\it
The  equalities
$$
W_{out} = \Theta^T W_{in}\Theta^T,
\quad S^{-1} = \Theta^T S\Theta^T
$$
are fulfilled.
}

\medskip
{\bf Corollary 2.4 ($PT$- and $CPT$--symmetry).}

{\it The equalities
$$
W_{out} = \Theta^T \Theta^P
W_{in}\Theta^T \Theta^P,
\quad S^{-1} =
 \Theta^T \Theta^P
S \Theta^T \Theta^P
$$
are fulfilled.
}

\medskip
{\bf Remarks.}

1. The symmetry $\Theta$ for the classical system means
the following.
$\Theta$ is a mapping of initial data,
$(\varphi,\pi)\to \Theta(\varphi,\pi),$ if $u$
is a solution of Eq. (2.1) with initial data
$(\varphi,\pi),$ then
$\Theta(\varphi,\pi)$ are
 initial data of some
solution $u^{\Theta}.$ In particular, if $u(t,x)$ is a
solution of Eq. (2.1), then $u^{\Theta^T}(t,x)=u(-t,x)$
is the corresponding solution of Eq. (2.1).

2. Since  $P$-symmetry
(the space reflection symmetry) is obvious, so
$T$-symmetry implies
$PT$-symmetry and $CPT$-symmetry for the nonlinear
equation  (2.1).
Of  course, $CPT$-symmetry is more fundamental.

3. On the other hand Theorem  2.2 is connected with
 the  correspondence between the time
reflection and the complex conjugation.
Theorem 2.2 appears in the fact that the interacting
quantum field is a Hermitian bilinear form. For the
bilinear form--solution this was pointed out in
\cite{Hei74},
\cite{Osi94a}.

4. The indexation of free solutions with the help of
positive frequency  part corresponds
 to the fact that the real diagonal of
Wick symbol of the interacting quantum field
 satisfies the classical nonlinear wave equation.

\medskip
{\it Proof of Theorem  2.3.
}

Let $u(t,x)$ be a solution of classical
nonlinear equation  (2.1) with finite
energy. Then
  $u(-t,x)$ (and $u(-t,-x)$) is a solution
of the same equation also.
 The
 solution $u$
defines the unique
 initial, $in$-, and
 $out$-data  at time zero,
\be
(\varphi(x), \pi(x))
&=&
(u(0,x), \dot u(0,x)),
\\
(\varphi_{in}(x), \pi_{in}(x))
&=&
 (u_{in}(0,x),\dot u_{in}(0,x)),
\\
(\varphi_{out}(x), \pi_{out}(x))
&=&
 (u_{out}(0,x), \dot u_{out}(0,x)),
\ee
where $u_{in}$ и $u_{out}$
are the unique free solutions such that
$$\Vert(u(t,\cdot),\dot u(t,\cdot)) - (u_{in}(t,\cdot),
\dot u_{in}(t,\cdot)) \Vert_{H^1\oplus L_2}
\quad\mbox{ for }\quad t\to-\infty$$
and
$$\Vert(u(t,\cdot),\dot u(t,\cdot))
- (u_{out}(t,\cdot),\dot
u_{out}(t,\cdot)) \Vert_{H^1\oplus L_2}
\quad\mbox{ for } \quad t\to+\infty.$$


This implies immediately
that the solution $u^{\Theta^T}$
 $(\, = u(-t,x)\,)$ has the initial data
equal to $(\varphi, -\pi)$
at time  $t=0,$
 the $in$--data equal to
$(\varphi_{out}, -\pi_{out}),$ and
 the $out$-data equal to
$(\varphi_{in},-\pi_{in}).$

Therefore, since by definition
 the wave operators satisfy the following equalities
$$
W_{in}(\varphi_{in},\pi_{in}) = (\varphi,\pi),
\quad
W_{out}(\varphi_{out},\pi_{out})
= (\varphi,\pi),
$$ $$
S(\varphi_{in},\pi_{in}) =
(\varphi_{out},\pi_{out}),
\quad S = W^{-1}_{out} W_{in},
$$
the previous assertion implies the equality
$$
\Theta^T W_{in}(\varphi_{in},\pi_{in})
= W_{out}\Theta^T(\varphi_{in},\pi_{in}), $$ i.e.
$$
W_{out} = \Theta^T W_{in}\Theta^T.
$$
The last equality implies that
$$
S^{-1}=\Theta^T S\Theta^T.
$$
Theorem 2.3 is proved.

\medskip
{\bf Remark.}

It is well known that $S$ is not the identity operator
\cite{MorS72b}.
This implies that
 $W_{out}\not=W_{in}.$



\section*{3. Wick kernel of the interacting
quantum field.}

The existence of classical wave operator and
 holomorphity
allow us to introduce the Wick kernel of
the interacting quantum field  and
to construct the bilinear form
given   by this kernel.

For this purpose we introduce the Hilbert Fock
space for the free quantum field $\phi_{in}.$
We introduce the Fock space
 as it described in
\cite{Rac75},
but instead of the basis used by
R\c{a}czka
\cite{Rac75}
we take a basis in the
complex Hilbert space
$H^{1/2}(\R^3,\C)$
 with explicit introduction
of its pure real and pure imaginary
parts  (for the integral by Paneitz, Pedersen,
Segal, Zhou
\cite{PanPSZ91},
see also
\cite{BaeSZ92},
we may use any basis in the complex
Hilbert space).

Let us take a
(canonical and standard) basis in
the space $H^{1/2}(\R^3,\C)$ analogous to the basis
introduced in
\cite[ch. 5.1, p. 141]{BogLT69}
(the space $H^{1/2}(\R^3,\C)$
corresponds to the dual space with dualization
$(.,.)_{L_2}$ of the one-\-particle
subspace of the Fock space,
this one-\-particle subspace is
 denoted by
Bogoliubov, Logunov, Todorov
\cite[ch. 5.1, p. 141]{BogLT69}
 as ${\cal H}_1$).

Similar to
\cite[ch.5.1, p. 141]{BogLT69}
we introduce our basis and note that all basic vectors
in our basis are
 pure real--valued
 functions in coordinate space.
Thus, we define for the (relativistic)
coordinate operator
(in momentum space)
$$ {\bf q}=\sqrt{\mu(p)}i\nabla_{{\bf p}}
{1\over\sqrt{\mu(p)}} =
i\nabla_{{\bf p}}-{i\over 2}{{\bf p}
\over m^2+{\bf p}^2}, $$
in this case the momentum operator
is the operator of multiplication on ${\bf p}$
(in momentum space).

The operators $q^j$ and $p^j,$ $j=1,2,3,$ satisfy
the canonical commutation relations
$$
[q^j,p^{j'}]=i\delta_{j,j'},\quad
[q^j,q^{j'}]=[p^j,p^{j'}]=0,$$
are hermitian (and essentially self--adjoint on
$\Sc(\R^3,\C)$)
in the Hilbert space  ${\cal H}_1$
with the scalar product
$$
(\psi_1,\psi_2)
= \int\overline{\psi_1(p)}\psi_2(p){d^3 p
\over\mu(p)},
$$
 the Hilbert space
  ${\cal H}_1$
coincides with
 the Hilbert space
  $H^{-1/2}(\R^3,\C)$ in coordinate space.

We define the operators
$$ b_j = {1\over\sqrt 2}(q^j+ip^j),
\quad b^\ast_j={1\over\sqrt 2}(q^j-ip^j),\quad
j=1,2,3.  $$
These operators are mutually adjoint and correspond
to the analogous operators in
\cite[ch. 5.1, p. 141]{BogLT69},
 $$ b_j=ib_j^{BLT}\;\;
\mbox{
(\cite[ch. 5.1, p.141]{BogLT69}),
}
$$ $$
b^\ast_j=-ib_j^{BLT}\;\;
\mbox{
(\cite[ch. 5.1, p. 141]{BogLT69}).
}
  $$
They satisfy
the  commutation relations
 $$
[b_j,b^\ast_{j'}]=\delta_{j,j'}, \quad
[b_j,b_{j'}]=[b^\ast_j,b^\ast_{j'}]=0.
$$
Let $e_0$ be the normalized vector,
satisfying the equalities
$$
b_j\; e_0=0,\quad j=1,2,3.
$$
In the considered momentum space
$$
e_0(p) =
 \pi^{-3/4}
(m^2+p^2)^{1/4} e^{-p^2/2}.
$$
The orthonormalized basic vectors in  ${\cal H}_1$
we define as
the Hermite functions (the Hermite polynomials)
defined on the mass hyperboloid,
$$
e_k \equiv e_{k_1 k_2 k_3} =  {(b^\ast_1)^{k_1}
(b^\ast_2)^{k_2}(b^\ast_3)^{k_3}\over
\sqrt{k_1! k_2! k_3!}} e_0,\quad k_1, k_2, k_3 = 0,1,...\, .
$$

In $x$--space the vector  $e_0(x)$ is
a real--valued function and the
vectors $e_k(x)$ are real--valued functions
also.
This follows from the fact, that
in the $p$--space
$$
b^\ast_j = {1\over\sqrt 2}
(\mu^{1/2}i\nabla_{p^j}\mu^{-1/2}-ip^j)=
{1\over\sqrt 2}
(\mu^{1/2}x^j\mu^{-1/2}-{\partial\over \partial x^j})
$$
and from the reality of operators $\mu$
as convolution operators.
Here
$$
q=\mu^{1/2}x \mu^{-1/2}
=\mu^{1/2}i\nabla_p\mu^{-1/2}
$$
is the operator of relativistic coordinate
in $p$--space (= in momentum space),
see
\cite[ch. 5.1, (2.5.6), p. 140]{BogLT69},
$\mu^{-1/2}=(-\Delta+m^2)^{-1/2}$ and gives
 an isomorphism  $L^2(\R^3) \to H^{1/2}(\R^3).$

In $x$-space the coordinate operator has the form
$$
q=\mu^{-1/2}x\mu^{1/2}
$$
because  the function $F(x)$
from the one-particle space of the Fock
space corresponds
to the function
 $$\mu({\cal F} F)(p)$$
 in the  $p$--space,
where  $\cal F$
is the Fourier transform. If, correspondingly,  $F(p)$
is a function that is square integrable over
$d^3 p/\mu(p),$
then in $x$--space this function corresponds
to the function
$\mu^{-1}({\cal F}F)(x)$
and the operator of relativistic coordinate is
 $\mu^{-1/2}x^j\mu^{1/2}$
(here all identifications appear  as
the consideration of the same solution
in different indexing spaces).

In this case, all functions that define the basis
are orthonormal and real--valued in  the $x$--space.
In the $p$--space these functions satisfy the
following relation
 $$ F(p)=\overline{F(-p)}.  $$

\medskip

The functions
$$
\varphi_k(x) =
c\int\exp(-i p x)
e_k(p)d\rho(p),
$$
$$
 c=[2(2\pi)^3]^{-1/2}, \quad
d\rho(p) = d^3 p/ \mu(p),
\quad \mu(p)=(p^2+m^2)^{1/2},
$$
are orthonormal real in the complex Hilbert
space
$ H^{1/2}(\R^3,\C)$ and generate a basis
in this Hilbert space.

We use the following
representation of the free (scalar, hermitian)
quantum field in the Fock space, see
\cite{BaeSZ92},
\cite[ch. 5.1]{BogLT69},
\cite{Rac75}.
We denote the Fock space for the free field as
$\Phi o\kappa$ and chose it as
$${\Phi o\kappa} =
\oplus^{\infty}_{n=0}{\cal H}_n,$$
where ${\cal H}_0=\C,$
 ${\cal H}_1=
 H^{-1/2}(\R^3,\C),$
 ${\cal H}_n=
\mbox{sym}\; \widehat\otimes_n H^{-1/2}(\R^3,\C).$
The Fock space
 $ {\Phi o\kappa}$
 is associated with the free quantum field
$\phi_{in}.$
Let
$\phi^{-}_{in}(0,x),$ $\phi^{+}_{in}(0,x)$
be the
annihilation and creation part of the free
quantum  field, i.e. the negative and
positive frequency part of the field $\phi_{in}$
at time zero. Using the chosen basis
$\{e_k(p)\}^\infty_{k=0}$ we define
the smeared  annihilation and creation operators
$a_k$ and $a^\ast_k,$
$$ a_k = \int \phi^{-}_{in}(0,x)
 e_k\widetilde{}(x) d^3x,
 \quad a^\ast_k = \int \phi^{+}_{in}(0,x)
 e_k\widetilde{}(x) d^3x,
$$
where
$ e_k\widetilde{}(x) $ is a basis vector in
coordinate space (i.e. the Fourier transform of the
vector
$ e_k(p) $ ).

The operators  $a_k,$ $a^\ast_k,$
 satisfy the canonical commutation relations
$[a_k,a^\ast_l]=\delta_{kl}$ and form
an infinite-dimensional
nilpotent Lie algebra, which
is irreducible in the Fock space
 $ {\Phi o\kappa}.$

The free quantum field $\phi_{in}(t,x)$ has the
following form in terms of annihilation $a_k$ and
 creation $a^\ast_k$ operators:
$$
\phi_{in}(t,x)
= \sum_k (\varphi_k(t,x)
a^\ast_k + \overline{\varphi_k(t,x)} a_k), $$
where
$$
\varphi_k(t,x) =
c\int\exp(i\mu(p) t-i p x)
e_k(p)d\rho(p),
$$

$$
d\rho(p) = d^3 p/\mu(p), \quad
\mu(p)=(p^2+m^2)^{1/2},\quad c=[2(2\pi)^3]^{-1/2},
$$
$$
\varphi_k(0,x)
 \equiv \varphi_k(x) = c\int\exp(-i p x)
e_k(p)d\rho(p),
$$ $$
\varphi_k(x)
= \overline{\varphi_k(x)},
$$
with our choice of the real basis,
$e_k(p) = \overline{e_k(-p)}.$
In this case
$$
\overline{\varphi_k(t,x)} = \varphi_k(-t,x).
$$

The vacuum averages of a sum of products of
the field $\phi_{in}$ coincide
with the corresponding value of the free
Wightman functional.

Now we construct  convenient dense subspaces
 contained in the Fock space. Let
$$
e_z \equiv e(z) = \exp(za^\ast)\Omega,
\eqno(3.1)
$$
$$
|z\rangle
= \exp(-{1\over 2}\Vert z\Vert^2)
\exp(za^\ast)\Omega ,
$$
where $$\Vert z\Vert^2 =
\sum_{k=0}^{\infty} |z_k|^2,\quad
za^\ast = \sum^\infty_{k=0} z_k
a^\ast_k, $$
$z$
is an element of the
complex Hilbert space  $H^{1/2}(\R^3,\C),$
$ z_k
= \langle \varphi_k, z\rangle_{H^{1/2}(\R^3,\C)}.
$
The vectors
$e_z,$
$z\in H^{1/2}(\R^3,\C),$
 are called
the coherent vectors and the vectors
$|z\rangle$ are called the
coherent state vectors.
It follows from
Eq. (3.1)
 that
$$ (\exp(z'
a^\ast)\Omega,\exp(za^\ast)\Omega)
=\exp(\langle z', z\rangle_{H^{1/2}(\R^3,\C)})
= \exp(\sum^\infty_{k=0} \bar z'_k z_k)
$$
and
\be \langle z'|z\rangle
&=&\exp\{-\sum^\infty_{k=0}
[{1\over2}|z'_k-z_k|^2-i\mbox{ Im }
(\bar z'_k z_k)]\}\cr
&=&\exp(\sum^\infty_{k=0}(\bar z'_k z_k-
{1\over 2}|z'_k|^2-{1\over 2} |z_k|^2))
\ee
in the chosen earlier basis
$$ z = \sum z_k\varphi_k, \quad z_k
= \langle \varphi_k, z\rangle_{H^{1/2}(\R^3,\C)},
$$ $$
\langle z_1, z_2\rangle_{H^{1/2}(\R^3,\C)}
= \sum^\infty_{k=0} \bar z_{1,k} z_{2,k}.
$$

Let $A$ be a subset in $H^{1/2}(\R^3,\C).$
Let $D_{coh}(A)$ be the subspace in the Fock
space of the $in$-field generated by finite linear
combinations of coherent vectors from $A,$
i.e. if $\chi \in D_{coh}(A),$ then
$$ D_{coh}(A)= \{\chi\in  D_{coh}(A)\, |\,
\chi=\sum \alpha_j e(z_j),\,\, \mbox{{\it the sum is
finite,}} \,\alpha_j\in\C, z_j\in A\}.
$$

We define the
Wick kernel for the interacting quantum field as
the  function
$(z_1,z_2)\to\phi(e_{z_1},e_{z_2})$ from
$ H^1(\R^3,\C)\times H^1(\R^3,\C)$ into
$ H^1(\R^3,\C),$
$$
\phi(e_{z_1},e_{z_2})=\exp(\langle z_1,
z_2\rangle_{H^{1/2}(\R^3,\C)})
{1\over 2}(RWR^{-1}(\bar z_1+z_2) +
\overline{RWR^{-1}(\overline{\bar z_1+z_2})}),
$$
where $RWR^{-1}$ maps   $H^1(\R^3,\C)$ into
$H^1(\R^3,\C).$
Thus,  the Wick kernel smoothed with some test function
over space is equal to
$$ \int(\phi(e_{z_1},e_{z_2}))(x)h(x)d^3 x =
\exp(\langle z_1,z_2\rangle)
$$
$$
{1\over 2}\Bigl(
\int(\mbox{ Re }
RWR^{-1}(\bar z_1+z_2))(x)h(x)d^3 x
 + i \int(\mbox{ Im }
RWR^{-1}(\bar z_1+z_2))(x)h(x)d^3 x $$
$$
+ \int(\mbox{ Re }
RWR^{-1}(\overline{\bar z_1+z_2}))(x)h(x)d^3 x
- i \int(\mbox{ Im }
RWR^{-1}(\overline{\bar z_1+z_2}))(x)h(x)d^3 x
\Bigr)
$$
$$ =\exp(\langle z_1,z_2\rangle)
{1\over 2}\Bigl(
\langle
\overline{RWR^{-1}(\bar z_1+z_2)},\mu^{-1}h
\rangle_{H^{1/2}(\R^3,\C)}
 +
\langle
RWR^{-1}(\overline{\bar z_1+z_2}),
\mu^{-1}h
\rangle_{H^{1/2}(\R^3,\C)}
\Bigr).$$
The smoothed Wick kernel corresponds
to the smoothed field
at the time zero,
i.e. to the bilinear form corresponding to
$\int\phi(0,x)h(x)d^3x.$

Analogously, we define the Wick kernel for the
$out$-going quantum field as
$$
\phi_{out}(e_{z_1},e_{z_2})=\exp(\langle z_1,
z_2\rangle_{H^{1/2}(\R^3,\C)})
{1\over 2}(RSR^{-1}(\bar z_1+z_2) +
\overline{RSR^{-1}(\overline{\bar z_1+z_2})}).
$$

In
\cite{Osi95b} 
we show that the interacting and $out$-going
Wick kernel can be reconstructed
 uniquely
with the help of the bilinear form--solution,
that has been constructed in
\cite{Osi94a}.
In this paper we
consider the introduced Wick kernels, their
properties  and the bilinear
form defined by these Wick kernels.

\medskip
{\bf Theorem 3.1.}

{\it The Wick kernel
$$
\phi(e_{z_1},e_{z_2})
=\exp(\langle z_1,z_2\rangle_{H^{1/2}(\R^3,\C)})
{1\over 2}(RWR^{-1}(\bar z_1+z_2)
+ \overline{RWR^{-1}\overline{(\bar z_1+z_2)}})
$$
is correctly defined on
$ H^1(\R^3,\C)$ $(\subset H^{1/2}(\R^3,\C))$
as the map from
$$ H^1(\R^3,\C)\times H^1(\R^3,\C)$$
into $H^1(\R^3,\C)$
and is complex
antiholomorphic on $z_1$ and
 complex
 holomorphic on $z_2.$
In particular, for $z_1, z_2$
belonging to finite--dimensional subspaces in
$H^1(\R^3,\C),$
$$z_1(\alpha_1)=\sum_{j=1}^n\alpha_{1,j}z_{1,j},
\quad z_2(\alpha_2)=\sum_{j=1}^n\alpha_{2,j}z_{2,j},
 \quad
z_{1,j}, z_{2,j} \in H^1(\R^3,\C),
\;\;\alpha_{1,j}, \alpha_{2,j}\in \C,
$$
$h\in H^1(\R^3,\C),$
the function
$$
\langle h,\phi (e(z_1(\alpha_1)),e(z_2(\alpha_2)))
\rangle_{H^{1/2}(\R^3,\C)}
$$
is an entire complex antiholomorphic function
on  $\alpha_1\in \C^n$ and
 an entire complex holomorphic function
on $\alpha_2\in \C^n$
 in the usual sense.
Furthermore,
$$
\overline{\langle h,
\phi (e_{z_1},e_{z_2})\rangle}_{H^{1/2}(\R^3,\C)}=
\langle \bar h,\phi (e_{z_2},e_{z_1})
\rangle_{H^{1/2}(\R^3,\C)},
$$
i.e.
$$
\overline{\phi (e_{z_1}, e_{z_2})}
= \phi (e_{z_2}, e_{z_1}),
$$
where a complex conjugation is defined as
the complex conjugation of a function
with values in $H^1(\R^3,\C).$
 It is valid the following estimate
$$
|\langle h,\phi(e_{z_1},
e_{z_2}
\rangle_{H^{1/2}(\R^3,\C)}
|\leq c
\Vert h\Vert_{L_2(\R^3,\C)}
\exp(\mbox{Re}(\langle z_1,z_2
\rangle_{H^{1/2}(\R^3,\C)}))
\;\Vert\bar z_1+z_2\Vert_{H^1(\R^3,\C)}.
$$

The analogous assertion is also valid for the $out$-going
field $\phi_{out}.$
}

\medskip
 The proof of Theorem 3.1
follows from the assertion of Theorem 2.1.

\medskip
{\bf Remark.}

The Hermitian symmetry is implied by
 the
explicit expression of the Wick kernel and
by  the equality
$$
W_{in} = \Theta^T W_{out}\Theta^T,
$$
or $$ \overline{RW_{in}R^{-1}z} =
R\Theta^T W_{in}\Theta^T R^{-1}\bar z
= R W_{out} R^{-1}\bar z
$$
for the wave operators
$W_{in},$ $W_{out}.$
 In particular, for any real test function  $h$
$$\int\phi(e_{z_1}, e_{z_2})(0,x)h(x)d^3 x =
\int\overline{\phi(e_{z_2}, e_{z_1})}(0,x)h(x)d^3 x,
$$
the integral is defined correctly
 (because
$
\phi(e_{z_1}, e_{z_2})\in H^1(\R^3,\C)
$
for $z_1, z_2\in H^1(\R^3,\C)$).

\medskip
This Wick kernel defines a bilinear form on
$D_{coh}(H^1(\R^3,\C))\times
D_{coh}(H^1(\R^3,\C)),$
i.e. on the subspace in the Fock space of the $in$--field
generated by finite linear combinations of
coherent vectors with
finite energy.

\medskip
{\bf Theorem 3.2.}

{\it Let  $\chi_1,\chi_2\in D_{coh}(H^1(\R^3,\C)),$
$\chi_1=\sum\alpha_{1,j}e_{z_{1,j}},$
 $\chi_2=\sum\alpha_{2,j}e_{z_{2,j}},$
then $$
\phi(\chi_1,\chi_2) =
\sum\overline{\alpha}_{1,j} \alpha_{2,j}
\phi(e_{z_{1,j}}, e_{z_{2,j}})
\eqno(3.2)
$$
is a bilinear form with values in $H^1(\R^3,\C).$
In addition, the expression
$$
\langle h,\phi(\chi_1,\chi_2)\rangle_{H^{1/2}(\R^3,\C)}
=\sum\overline{\alpha}_{1,j}\alpha_{2,j}
\langle h,\phi(e_{z_{1,j}},
e_{z_{2,j}})\rangle_{H^{1/2}(\R^3,\C)}
$$
defines the bilinear form also.
The bilinear forms
$\phi(\chi_1,\chi_2)$ and
$\langle h,\phi(\chi_1,\chi_2)\rangle_{H^{1/2}(\R^3,\C)}$
 are Hermitian symmetric
$$\overline{\phi(\chi_1,\chi_2)} =
\phi(\chi_1,\chi_2),
$$
$$\overline{\langle h,
\phi(\chi_1,\chi_2)\rangle}_{H^{1/2}(\R^3,\C)}
= \langle \overline{h},
\phi(\chi_1,\chi_2)\rangle_{H^{1/2}(\R^3,\C)}.
$$

The analogous assertion is also valid for the
$out$-going quantum field $\phi_{out}.$
}

\medskip
\medskip

Let $H_{in}$ be the Hamiltonian of the (free quantum)
$in$-coming field $\phi_{in}.$ Due to the
energy conservation it is equal to the total
Hamiltonian.

\medskip
{\bf Theorem 3.3
(Energy positivity and holomorphity).}

{\it Let  $\chi_1=\sum\alpha_{1,j}e_{z_{1,j}},$
 $\chi_2=\sum\alpha_{2,j}e_{z_{2,j}},$
$\chi_1,\chi_2\in D_{coh}(H^1(\R^3,\C)).$
The expression
 $$
\phi(\exp(it_1H_{in} - s_1H_{in})\chi_1,
 \exp(it_2H_{in} - s_2H_{in})\chi_2)
$$ $$ =
\sum\overline{\alpha}_{1,j}\alpha_{2,j}
\phi(e(\exp(i\mu t_1 - \mu s_1)z_{1,j}),
 e(\exp(i\mu t_2 - \mu s_2)z_{2,j}))
\eqno(3.3)
$$
is correctly defined for
$\mbox{ Im }s_1\geq 0,\mbox{ Im }s_2\geq 0$
as a bilinear form on
 $$D_{coh}(H^1(\R^3,\C))\times D_{coh}(H^1(\R^3,\C))$$
with values in
$H^1(\R^3,\C).$
Expression (3.3) depends
 antiholomorphically on
$t_1+is_1$ and
 holomorphically
on $t_2+is_2$).

The analogous assertion is also valid for the
$out$-going quantum field $\phi_{out}.$
}

\medskip
\medskip
{\it Proof of Theorem 3.2 and 3.3.}
It is evident that Theorem 3.2 is the consequence
of Theorem 3.3 with
$t_1 = t_2 = s_1 = s_2 = 0.$
Therefore, it is sufficient to prove Theorem 3.3.

The assertion
about  holomorphity on
$\overline{t_1+is_1},$ $t_2+is_2$
 is implied by
 Theorem 3.1 (or by Theorem 2.1) and
by  holomorphity of
the operator
$\exp(i\mu t - \mu s)=
\exp i\mu(t+is)$ on $t+is$ for $s>0.$

\medskip
Now we prove that the expression
$\phi$ is a correctly defined bilinear form.
We consider here two variants of the proof for the
interacting field $\phi.$ The field $\phi_{out}$
can be considered analogously.

\medskip
The first variant of the proof uses
 Theorem 1
\cite[Theorem 1]{PanPSZ91}.
The Wick kernel is given by
nonlinear operators
$RWR^{-1}$ and $R\Theta^T W\Theta^T R^{-1}.$
 The complex holomorphity
of the considered Wick kernel was proved in
\cite{Osi95a}.
In order to use Theorem 1
\cite[Theorem 1]{PanPSZ91}
we take as an
operator $B$ the operator
 $$
B=\mu+\mu^{-1/2}(1+x^2)\mu^{1/2}.
$$
 It is evident, that this operator is
positive in the coordinate space
$H^{1/2}(\R^3,\C),$
self-adjoint, and $\exp(-sB)$
is a  trace class operator for  $s>0.$

\medskip
{\bf Remark.}

 It is possible to use more simple choice
of (positive, self-adjoint) operator
 $B,$ for instance,  $$
B=\mu^\alpha+\mu^{-1/2}(1+x^2)^\beta\mu^{1/2},
$$
or
$$
B=\mu^{\alpha-1/2}(1+x^2)^\beta\mu^{\alpha+1/2},
$$
for sufficiently large
$\alpha,\beta,$ such that the operator $B^{-1}$
 is  nuclear (as a product of two
Hil\-bert-\-Schmidt operators).

\medskip
We now consider the Wick kernel
$\langle h,\phi(e(e^{-s_1 B}z_1),
e(e^{-s_2 B}z_2)\rangle.$
It is easily
to see that this kernel
satisfies the estimate (2) of Theorem 1
\cite[Theorem 1]{PanPSZ91}.
This is implied by the explicit form of
the Wick kernel,  conservation of energy
(which gives an estimate for the nonlinear
terms), and by a simple
estimate
$B\leq c(s)\exp(sB),$ $s>0.$
Therefore, the conditions of
Theorem 1
\cite[Theorem 1]{PanPSZ91}
are fulfilled and $\phi$ is a correctly
defined bilinear form.
It is possible to extend
this bilinear form by continuity
on the subspace generated by finite
linear combinations of coherent vectors
with finite energy.

The first proof of Theorem 3.3 and
Theorem 3.2 is complete.

\medskip
\medskip
{\it The second  proof of Theorem 3.3.}
 Let $\chi_1(z),$ $\chi_2(z)$
 be the holomorphic functions on the
space $H^1(\R^3,\C)$  ($\subset H^{1/2}(\R^3,\C)$)
corresponding to the vectors
  $\chi_1,$ $\chi_2$
in the complex wave representation of the
Fock space
(see
\cite[ch. 1, Theorem 1.13, p. 67-68]{BaeSZ92}.
 Let $\theta$ be strictly positive, $\theta > 0,$
 and such that for $\beta_1, \beta_2 \in \C,$
 $ |\beta_1| + |\beta_2| < \theta,$
the Taylor series
 in
$\beta_1, \beta_2 $ for the functions
$\phi(e(\beta_1 z_{1,j_1}),
e(\beta_2 z_{2,j_2}))$
 converge in $H^1(\R^3,\C)$
for all $j_1, j_2$
($\chi_1,$ $\chi_2$ are finite linear combinations
of coherent vectors).
Let $$\chi_1(\beta_1)=\chi_1(\beta_1z), \quad
\chi_2(\beta_2)=\chi_2(\beta_2z).
$$
For our sufficiently small
$\beta_1, \beta_2,$
 $ |\beta_1| + |\beta_2| < \theta,$
we write the following equality
 $$
 \sum_{j_1,j_2}
\overline{\alpha}_{1,j_1}\alpha_{2,j_2}
\phi(e(\beta_1 z_{1,j_1}),
 e(\beta_2 z_{2,j_2}))
$$
 $$
= \sum_{j_1,j_2}
\overline{\alpha}_{1,j_1}\alpha_{2,j_2}
\exp(\langle\beta_1 z_{1,j_1},\beta_2 z_{2,j_2}\rangle)
\sum_{n=1}^{\infty}
\phi_n(
\underbrace{
(\overline{\beta}_1\overline{z}_{1,j_1} +
\beta_2 z_{2,j_2})
\otimes...\otimes
(\overline{\beta}_1\overline{z}_{1,j_1}
+ \beta_2 z_{2,j_2})}_{n})
$$
$$
= \lim_{N\to\infty}\lim_{V,\sigma}
\sum_{n=1}^N \sum_{j_1,j_2}
\overline{\alpha}_{1,j_1}\alpha_{2,j_2}
\exp(\langle\beta_1 z_{1,j_1},\beta_2 z_{2,j_2}\rangle)
$$ $$\phi_n(
\underbrace{
(\overline{\beta}_1\overline{z}_{1,j_1,V,\sigma}
+ \beta_2 z_{2,j_2,V,\sigma})
\otimes...\otimes
(\overline{\beta}_1\overline{z}_{1,j_1,V,\sigma}
+ \beta_2 z_{2,j_2,V,\sigma})}_{n}).
\eqno(3.4)
$$
Here $\sigma$ is an ultraviolet cut-off and $V$ is a space
cut-off (i.e. $\sigma$ is the smoothing with some
test function and $V$ is the multiplication on some test
function, $\sigma$ tends to the $\delta$-function and
$V$ tends to 1). Further,
$$\phi_n = {1\over 2}(R_n + \overline{R}_n)
$$
is a tempered distribution (i.e. it
is a generalized function from the Schwartz space),
$R_n$ is the tempered distribution defined
uniquely by the $n$-linear continuous form
$$
{1\over n!}d^n RWR^{-1}(0)
$$
 on
$H^1(\R^3,\C) \otimes... \otimes
 H^1(\R^3,\C)$
and by the  Schwartz
nuclear theorem.
We denote this
 $n$-linear form and the generalized
function by the same notation $R_n.$
The generalized function $R_n$
is such that
$$ R_n(f_1\otimes...\otimes f_n) =
{1\over n!}d^n RWR^{-1}(0)(f_1,...,f_n).
$$

A Wick polynomial is a correctly defined bilinear form,
for instance, on
$$D_{coh}(\Sc(\R^3,\C)) \times D_{coh}(\Sc(\R^3,\C)).$$
This is implied easily by the explicit form of Wick
 monomial
$$
:\phi_{in}(f_1)
...\phi_{in}(f_n):(\chi_1,\chi_2)=
2^{-n/2}\sum_{K\subseteq\{1,...,n\}}
\langle\prod_{k\in K}a(f_k)\chi_1,
\prod_{k\in\{1,...,n\}\backslash K}
a(f_k)\chi_2\rangle,
$$
 where $a$ is an annihilation operator.
This explicit form allows to extend
these expressions
by
continuity and linearity
 on tempered generalized functions,
see, for instance,
Reed, Simon
\cite[v. 2]{ReeS75},
\cite[Theorem 3 with appropriate operator $A$]{Bae89}.

Therefore, the expressions
$$
\phi_n(:\phi_{in}...\phi_{in}:),
$$
$\phi_n\in \Sc^{\prime}(\R^{3n}, H^1(\R^3,\C)),$
is  correctly defined as a bilinear form on
$$D_{coh}(\Sc(\R^3,\C)) \times D_{coh}(\Sc(\R^3,\C))$$
with values in  $H^1(\R^3,\C).$
On vectors
$\chi_1, \chi_2 \in
 D_{coh}(\Sc(\R^3,\C))$
this bilinear form is equal
 to
$$
\phi_n(:\phi_{in}...\phi_{in}:)
(\chi_1, \chi_2)
$$ $$
=
\sum_{j_1,j_2}\overline{\alpha}_{1,j_1}
\alpha_{2,j_2}
\exp(\langle z_{1,j_1}, z_{2, j_2}
\rangle_{H^{1/2}(\R^3,\C)})
 \phi_n((\overline{z}_{1,j_1} + z_{2,j_2})
\otimes...\otimes
(\overline{z}_{1,j_1} + z_{2,j_2})).
$$
Now if
  $\chi_1, \chi_1^\prime,
  \chi_2, \chi_2^\prime,
   \in
 D_{coh}(\Sc(\R^3,\C)),$
  $$\chi_1 = \chi_1^\prime, \quad
  \chi_2 = \chi_2^\prime, $$
we define
  $\chi_1(\beta_1), \chi_1^\prime(\beta_1),
  \chi_2(\beta_2), \chi_2^\prime(\beta_2).$
 It is clear that
  $$\chi_1(\beta_1) = \chi_1^\prime(\beta_1),
  \quad \chi_2(\beta_2) = \chi_2^\prime(\beta_2).$$
  Thus, for sufficiently small $\theta,$
 $|\beta_1| + |\beta_2| < \theta,$
Eq. (3.4) implies that
$$
\phi(\chi_1(\beta_1),\chi_2(\beta_2))
= \lim_{N}\lim_{V,\sigma}
\sum_{n=1}^N
\phi_n(
(\underbrace{
:\phi_{in}...\phi_{in}:}_{n})
(\chi_1(\beta_1),\chi_2(\beta_2))
$$
$$
=
 \lim_{N}\lim_{V,\sigma}
\sum_{n=1}^N
\phi_n(
(\underbrace{
:\phi_{in}...\phi_{in}:}_{n})
(\chi_1^\prime(\beta_1),\chi_2^\prime(\beta_2))
 = \phi(\chi_1^\prime(\beta_1),\chi_2^\prime(\beta_2))
$$
i.e.
$$
\phi(\chi_1(\beta_1),\chi_2(\beta_2)) =
\phi(\chi_1^\prime(\beta_1),\chi_2^\prime(\beta_2))
\eqno(3.5)
$$
for small $\beta_1, \beta_2.$
Since
$
\phi(\chi_1(\beta_1),\chi_2(\beta_2))
$ and
$\phi(\chi_1^\prime(\beta_1),\chi_2^\prime(\beta_2))
$
are entire holomorphic functions on
$\bar{\beta}_1,$
${\beta}_2$
Eq. (3.5) for small
$\beta_1,$ $\beta_2$ implies that it is
fulfilled for all
$\beta_1, \beta_2 \in \C, $ in particular, for
$\beta_1 =  \beta_2 = 1, $
that is,
$$
\phi(\chi_1,\chi_2)
=\phi(\chi_1(1),\chi_2(1))
= \phi(\chi_1^\prime(1),\chi_2^\prime(1))
= \phi(\chi_1^\prime,\chi_2^\prime).
$$

The proof of linearity in
$\chi_1, \chi_2 $
is analogous.

Therefore, (3.2) and (3.3) are
correctly defined bilinear forms.
The second proof of Theorem  3.2 and Theorem 3.3 is complete.

\medskip
{\bf Theorem 3.4
(Poincar{\'e} invariance).}

{\it The Wick kernel satisfies the following equalities
$$
\phi(e_{z_1},e_{z_2})((a,\Lambda)^{-1}(t,x))=
\phi(e(RU_0(a,\Lambda)R^{-1}z_1),
e(RU_0(a,\Lambda)R^{-1}z_2))(t,x),
$$
$$
\phi_{out}(e_{z_1},e_{z_2})((a,\Lambda)^{-1}(t,x))=
\phi_{out}(e(RU_0(a,\Lambda)R^{-1}z_1),
 e(RU_0(a,\Lambda)R^{-1}z_2))(t,x).
$$
Here $(a,\Lambda)\in {\cal P}_0,$
$$\phi(e_{z_1},e_{z_2})((t,x)) =
\phi(e(RU_0(-t)R^{-1}z_1),
e(RU_0(-t)R^{-1}z_2))(0,x),
$$ $$\phi_{out}(e_{z_1},e_{z_2})((t,x)) =
\phi_{out}(e(RU_0(-t)R^{-1}z_1),
e(RU_0(-t)R^{-1}z_2))(0,x),
$$
$ U_0(t) = U_0(a,\Lambda)|_{(a,\Lambda)=((t,0),1)}.
$
}

\medskip
\medskip
{\it Proof of Theorem 3.4.}
 Theorem 3.4 is the direct consequence of Theorem 2.2.
Theorem 3.4 is proved.


\section*{4. Conclusion.}

In this paper with the help of the complex
structure of the classical
nonlinear wave equation
we introduce the Wick kernel of the interacting
quantum field. In the next paper
\cite{Osi96b}
we prove that the constructed Wick kernel defines the
operator-valued generalized function from
  a Gelfand space
$\Sc^{\alpha}(\R^4),$ $\alpha>1$
\cite{GelS58}.
The proof uses mainly the technique of works
mentioned in the references.

\medskip
\medskip
\medskip
\medskip
\section*{Acknowledgments.}
\medskip
\medskip

This is the second paper of the project
 $\phi^4_4 \cap M.$
This paper is supported partly by RBRF 96-01649.

We acknowledge
Anatoly Kopilov, Valery Serbo,
 Vasily Serebryakov,
 Ludwig Faddeev, Anatoly Vershik,
 Peter Osipov,  Volja Heifets,  Julia
 and Zinaida
 for the help, advice, and criticism.


\begin{thebibliography}{999}

 \bibitem{Aiz82}Aizenman, M.: Geometric analysis
of $\phi^4$ fields and Ising models
(Part I \& II). Commun. Math. Phys. {\bf 86} (1982) 1-48.

 \bibitem{AizG83}Aizenman, M., Graham, R.: On the renormalized
coupling constant and the susceptibility in $\phi_4^4$
field theory and the Ising model in four dimensions.  Nucl.
Phys. {\bf B225 [FS9]} (1983) 261-288.

\bibitem{Bae89}
Baez, J.C.:
Wick  products of the free Bose field.
J. Funct. Anal. {\bf 86} (1989) 211-225.

\bibitem{Bae92}Baez, J.:
Scat\-ter\-ing and comp\-le\-te in\-teg\-r\-abi\-lity
in four di\-men\-si\-ons.
Mathematical aspects of classical field theory.
Proceedings of the
AMS-\-IMS-\-SIAM. Joint Sum\-mer Re\-se\-arch
Con\-fe\-ren\-ce held
July 20-\-26, 1991.
Ed. M.J.Go\-tay, J.E.Mar\-s\-den, V.Mon\-cri\-ef.
In: Contemporary Mathematics, {\bf v. 132}.
Rhode Is\-land, AMS Pro\-vi\-den\-ce, 1992,
pp. 99-116.

\bibitem{BaeSZ92}
Baez, J.,  Segal, I.E., Zhou, Z.:
Introduction to algebraic quantum field theory.
Princeton University Press, 1992.

 \bibitem{Ber65}
Berezin, F.A.:
The method of second quantization.
New York, Academic Press, 1966.

\bibitem{BogLT69}
 Bogoliubov, N.N., Logunov, A.A., Todorov, I.T.:
The foundations of axiomatic approach to quantum
 field theory.
Moscow, Nauka, 1969.


\bibitem{BogLOT87}
 Bogoliubov, N.N., Logunov, A.A., Oksak, A.I.,
Todorov, I.T.:
The basic principles of quantum field theory.
Moscow, Nauka, 1987.

\bibitem{BogMP58}
 Bogoliubov, N.N., Medvedev,~B.V., Polivanov,~M.K.:
Voprosi teorii dispersionnih sootnoshenij.
Moscow, Fizmatgiz, 1958.

\bibitem{BoSh73}
Bogoliubov, N.N., Shirkov, D.V.:
Introduction to the theory of quantized fields.
New York, Interscience, 1959.

 \bibitem{Cal88}Callaway, D.:
Triviality pursuit: can
elementary scalar particles exist?
Phys. Reports {\bf 167} No 5 (1988) 241-320.

\bibitem{ConT79}
Constantinescu, F., Thalhaimer, W.:
Ultradistributions and quantum fields:
Fourier-Laplace transforms and
boundary values of analytic functions.
Rept. on Math.Phys.
vol.16, No 2 (1979) 167-180


 \bibitem{Fro82}Fr\"ohlich, J.: On the triviality of
$\lambda\varphi_d^4$ theories and the approach
to the critical point in $d\ge 4$ dimensions.
Nucl. Phys. {\bf B200 [FS4]}
(1982) 281-296.

\bibitem{GelS58}
Gelfand I.M., Shilov G.E.:
Prostranstva osnovnikh i obobshchennikh funktsij
(Obobshchennie funktsii. Vipusk 2).
Moskwa, Fizmatgiz, 1958.

\bibitem{GelV61}
Gelfand, I.M., Vilenkin, N. Ya.:
Nekotorie primeneniya garmonicheskogo analiza.
Osnashchennie gilbertovi prostranstva.
(Obobshchennie funktsii. Vipusk 4).
Moskwa, Fizmatgiz, 1961.

\bibitem{Goo64}
Goodman, R. W.:
One-sided invariant subspaces and domains of
uniqueness for hyperbolic equations.
Proc.  Amer.  Math.  Soc.  {\bf 15} (1964) 653-660.

 \bibitem{Hei74}
Heifets, E.P.:
The $\phi^4_4$  classical wave equation
and the construction of the quantum field
as a bilinear form in the Fock space.
Novosibirsk, Institute of Mathematics, 1974.

\bibitem{Jos67}
Jost, R.:
The general theory of quantized fields.
Providence, Rhode Island, American Math. Soc., 1965.


\bibitem{MorS72a}
Morawetz, C. S. and Strauss, W. A.:
Decay and  scattering of solutions of
a nonlinear relativistic wave equation,
 Comm. Pure Appl. Math.  {\bf 25} (1972) 1-31.

\bibitem{MorS72b}
Morawetz, C. S. and Strauss, W.A.:
On a nonlinear scattering operator,
 Comm. Pure Appl. Math. {\bf 26} (1972) 47-54.


\bibitem{Osi84}
Osipov, E.P.:
Complex structure of the manifold of solutions
of classical wave equation and quantization.
6-th International Symposium on Information
Theory, Part III. Moscow-Tashkent, 1984,
pp. 164-166.

 \bibitem{Osi94a}
Osipov, E.P.:
Quantum interaction $\phi^4_4$: the construction of
quan\-t\-um field defined as a bilinear form.
Novosibirsk, Institute of Mathematics, TPh-205, 1994,
hep-th/9602003.

\bibitem{Osi94b}
Osipov, E.P.:
Quantum interaction $\phi^4_4$:
the construction of Wightman functions.
Novosibirsk, Institute of Mathematics, TPh-206, 1994.

\bibitem{Osi94MIT}
Osipov, E.P.:
Quantum interaction $\phi^4_4$:
the construction of the solution of quant\-um
wave equation,
the construction of Wightman functions.
Symposium on quan\-ti\-za\-ti\-on and
nonlinear wave equations,
MIT, Cambridge, June 1994.

\bibitem{Osi95a}
Osipov, E.P.:
Complex structure and solutions of
classical non-linear equation with
the interaction $u^4_4.$
Novosibirsk, Institute of Mathematics, TPh-207, 1995,
funct-an/9602002, submitted to J. Funct. Anal.

 \bibitem{Osi95b}
Osipov E.P.:
Quantum interaction $:\!\phi^4_4\!:.$
The construction of the quantum field as the solution
of the nonlinear quantum wave equation.
Novosibirsk, Institute of Mathematics,  TPh-210, 1995.


\bibitem{Osi96a}
Osipov, E.P.:
Complex structure and solutions of
the non-linear Klein-Gordon equation with
the interaction $u^4_4.$
Novosibirsk, Institute of Mathematics, TPh-212, 1996.
Submitted to Electronic Research Announcements of
the AMS.

\bibitem{Osi96b}
Osipov, E.P.:
The $\!:\phi^4_4:\!$ quantum field theory, II.
Integrability of Wick kernels.
In preparation.

\bibitem{Osi96c}
Osipov, E.P.:
In preparation.

\bibitem{Pan80-82}
Paneitz, S.:
Hermitian structures on solution varieties
of nonlinear relativistic wave equations,
{\it in:} ``Differential geometric methods
in mathematical physics (Clausthal, 1980)", pp. 108-118,
Lecture Notes in Math. 905, Springer, Berlin-\-New York,
1982.

\bibitem{Pan81}
Paneitz, S.:
Unitarization of symplectics and differential equations
in Hilbert space.
 J.  Funct.  Anal.  {\bf 41} (1981) 315-326.

\bibitem{Pan82}
Paneitz, S.:
Essential unitarization of symplectics
and applications to field quantization,
 J. Funct. Anal. {\bf 48} (1982) 310-359.

\bibitem{PanS80}
Paneitz, S.M. and Segal, I.E.:
Quantization of wave equations and hermitian structures
in partial differential varieties,
 Proc. Nat. Acad. Sci. USA {\bf N12} (1980) 6943-6947.

\bibitem{PanPSZ91}
Paneitz, S.M., Pedersen, J., Segal, I.E., Zhou, Z.:
Singular operators on boson fields as
forms on spaces of entire functions  on Hilbert space.
J. Funct. Anal. {\bf 100} (1991) 36-58.

\bibitem{PedSZ92}
Pedersen, J., Segal, E., Zhou, Z.:
Massless
$\phi^q_d$  quantum field
theories and the non\-tri\-vi\-al\-ity
of $\phi^4_4.$ Nucl. Phys.
{\bf B376} (1992) 129-142.

 \bibitem{Rac75}
R\c{a}czka, R.:
The  construction of solution of
nonlinear relativistic
wave equation in $\lambda:\!\phi^4\!:$ theory.
J. Math. Phys. {\bf 16}
(1975) 173-176.

\bibitem{RacS79}
R\c{a}czka, R., Strauss, W.:
Analyticity properties of
the scattering
operator in nonli\-ne\-ar relativistic classical
and prequantized
field theories. Rept. Math. Phys.
{\bf 16} No 3 (1979) 317-327.

\bibitem{ReeS72}
Reed, M., Simon, B.:
Methods of  modern mathematical physics.
Vol. I. Functional analysis.
New York, London, Academic Press, 1972.

\bibitem{ReeS75}
Reed M., Simon, B.:
Methods of Modern Mathematical Physics.
Vol. II. Fourier Analysis, Self-Adjointness.
New York, Academic Press, 1975.

\bibitem{ReeS79}
Reed, M., Simon, B.:
Methods of  modern mathematical physics
Vol. III.  Scattering theory.
New York, San Francisco, London,
Academic Press, 1979.

\bibitem{Seg66}
Segal, I.E.:
Dispersion for non-linear relativistic equations.
{\it In:} ``Proc. Conf. Math. Th. El. Particles",
 pp.79-108,
M.I.T. Press, Cambridge, MA, 1966.

\bibitem{Seg68}
Segal, I.E.:
Dispersion for non-linear relativistic equations, II.
 Ann. Sci. {\' E}cole Norm. Sup.
{\bf 1} (4) (1968) 459-497.

\bibitem{StrW64}
Streater, R.F., Wightman, A.S.:
PCT, spin and statistics, and
all that. New York-Amsterdam, Benjamin 1964.

\bibitem{Vla64}
Vladimirov, V.S.: Methods of the
theory of functions of many
complex variables. Cambridge, MIT Press 1966.


\end{thebibliography}
\end{document}